\documentclass[1p]{elsarticle}
\usepackage{amssymb}
\usepackage{amsfonts}

\newtheorem{theorem}{Theorem}
\newtheorem{conjecture}[theorem]{Conjecture}
\newtheorem{corollary}[theorem]{Corollary}

\newtheorem{observation}[theorem]{Observation}
\newtheorem{claim}
{Claim}

\newproof{pf}{Proof}

\begin{document}
\title{Distant set distinguishing total colourings of graphs}

\author{Jakub Przyby{\l}o\fnref{fn1,fn2}}
\ead{jakubprz@agh.edu.pl, phone: 048-12-617-46-38,  fax: 048-12-617-31-65}

\fntext[fn1]{Financed within the program of the Polish Minister of Science and Higher Education
named ``Iuventus Plus'' in years 2015-2017, project no. IP2014 038873.}
\fntext[fn2]{Partly supported by the Polish Ministry of Science and Higher Education.}

\address{AGH University of Science and Technology, al. A. Mickiewicza 30, 30-059 Krakow, Poland}

\begin{abstract}
The Total Colouring Conjecture suggests that $\Delta+3$ colours ought to suffice in
order to provide a proper total colouring of every graph $G$ with maximum degree $\Delta$.
Thus far this has been confirmed up to an additive constant factor,
and the same holds even
if one additionally requires
every pair of neighbours in $G$ to differ with respect to the sets of their incident colours,
so called \emph{pallets}.
Within this paper we conjecture that an upper bound of the form $\Delta+const.$ still remains valid even after extending the distinction
requirement to pallets associated with 
vertices at distance at most $r$, 
if only $G$ has minimum degree $\delta$ larger than a constant dependent on $r$.
We prove that such assumption on $\delta$ is then unavoidable and exploit the probabilistic method in order to
provide two supporting results for the conjecture. 
Namely, we prove the upper bound $(1+o(1))\Delta$ for every $r$,
and show that the conjecture holds if $\delta\geq \varepsilon\Delta$ for any fixed $\epsilon\in(0,1]$ and $r$,
i.e., in particular for regular graphs.
%
\end{abstract}

\begin{keyword}
Zhang's Conjecture \sep adjacent vertex distinguishing total chromatic number \sep total neighbour distinguishing index \sep $d$-strong total chromatic number \sep $r$-adjacent strong total chromatic number \sep $r$-distant set distinguishing total number \sep total neighbour distinguishing index by sums \sep Total Colouring Conjecture
\end{keyword}

\maketitle

\section{Introduction}
In~\cite{Zhang} the following intriguing extension of proper edge colourings was introduced.
Given a graph $G=(V,E)$, where by $E(v)$ we shall understand the set of edges incident with a vertex $v$ in $G$,
and its edge colouring $c:E\to\{1,2,\ldots,k\}$,
denote by
\begin{equation}\label{edgeScv}
S_c(v):=\{c(e):e\in E(v)\}
\end{equation}
the \emph{pallet of colours} incident with $v\in V$.
This shall also be denoted by $S(v)$ if $c$ is unambiguous from context.
A proper edge colouring $c$ of $G$ is called \emph{neighbour set distinguishing} 
or \emph{adjacent strong} 
if 
$S(u)\neq S(v)$ for every edge $uv\in E$.
The least number of colours in such a colouring is called
the \emph{neighbour set distinguishing index} or the \emph{adjacent strong chromatic
index} and denoted by $\chi'_a(G)$,
see~\cite{Akbari,BalGLS,FlandrinMPSW,Hatami,Zhang},
also for other notations used.
Surprisingly, it was conjectured in \cite{Zhang} that just one more colour than stemming from the
Vizing's Theorem on the sufficient number of colours to assure a proper edge colouring of a graph 
is (almost) always enough to distinguish neighbours by colour pallets as well.
\begin{conjecture}[\cite{Zhang}]\label{ZhangsConjecture}
For every connected graph $G$, $\chi'_{a}(G)\leq \Delta(G)+2$, unless $G$ is isomorphic to $K_2$ or $C_5$.
\end{conjecture}
The best general result corresponding to this conjecture is thus far the result of Hatami:
\begin{theorem}[\cite{Hatami}]
If $G$ is a graph with no isolated edges and maximum degree $\Delta>10^{20}$, then $\chi'_{a}(G)\leq \Delta+300$.
\end{theorem}
Moreover, $\chi'_a(G)\leq 3\Delta(G)$ by \cite{Akbari}, and $\chi'_a(G)\leq \Delta(G)+O(\log\chi(G))$ by \cite{BalGLS}.
Conjecture~\ref{ZhangsConjecture} was also verified e.g. for  bipartite graphs and for graphs of maximum degree 3, see~\cite{BalGLS}.

We shall focus on a correspondent of the concept above,
but 
embedded
in total colourings environment.
%
%
In this case, i.e., when $c:E\cup V\to\{1,2,\ldots,k\}$ is a proper total colouring
(i.e., no two adjacent vertices get the same colour,
no two incident edges get the same colour, and no edge gets the same colour as one of its endpoints),
the colour pallet of $v$, $S_c(v)$,
shall be understood slightly differently for every $v\in V$, namely
\begin{equation}\label{totalScv}
S_c(v):=\{c(e):e\in E(v)\}\cup\{c(v)\}
\end{equation}
then.
The total colouring $c$ is called \emph{adjacent vertex distinguishing}
if $S_c(u)\neq S_c(v)$ for every edge $uv\in E$.
The least number of colours in such a colouring is called
the \emph{adjacent vertex distinguishing total chromatic number} 
and denoted by $\chi''_a(G)$.
It was first considered in~\cite{Zhang_total}.
\begin{conjecture}[\cite{Zhang_total}]\label{ZhangsConjecture_total}
For every graph $G$, $\chi''_{a}(G)\leq \Delta(G)+3$.
\end{conjecture}
Note that again 
the conjectured upper bound for this parameter exceeds only by one
the expected corresponding upper bound for the necessary number of colours
in a proper total colouring of a graph $G$, $\chi''(G)$. 
\begin{conjecture}[The Total Colouring Conjecture]
For every graph $G$, $\chi''(G)$ $\leq$ $\Delta(G)$ $+2$.
\end{conjecture}
Answering the question formulated in Conjecture~\ref{ZhangsConjecture_total} above seems however even
more challenging than verifying its correspondent concerning just edge colourings,
as already the Total Colouring Conjecture, posed 
by Vizing~\cite{Vizing2} in 1968
and independently by Behzad~\cite{Behzad} in 1965 is still far from being fully settled
(unlike the edge colouring case, solved completely by Vizing).
The best result concerning this was delivered in 1998 by Molloy and Reed \cite{MolloyReedTotal}, who designed a complex probabilistic argument implying that $\chi''(G)\leq\Delta(G)+const.$.
\begin{theorem}[\cite{MolloyReedTotal}]\label{MolloyReedTh}
There exist constants $\Delta_0$ and $C_1$ such that for every graph $G$ with $\Delta(G)\geq\Delta_0$, $\chi''(G)\leq\Delta(G)+C_1$.
In particular, $C_1\leq 10^{26}$.
\end{theorem}
%
%
%
Basing on this result and the approach of Hatami~\cite{Hatami},
quite recently Coker and Johanson proved the following.
\begin{theorem}[\cite{CokerJohanson}]
There exists a constant $C'$ such that for every graph $G$, $\chi''_a(G)\leq\Delta(G)+C'$.
\end{theorem}
It is also known that Conjecture~\ref{ZhangsConjecture_total} holds e.g. for
cycles, complete graphs, complete bipartite graphs and trees, see~\cite{Zhang_total},
graphs of maximum degree at most three, see e.g~\cite{Wang_totalDelta3},
$K_4$-minor free graphs, see~\cite{WangWangK_4_total_set}, planar graphs with $\Delta(G)\geq 14$~\cite{WangHuang_planar}, and for graphs with `small' average degree, see~\cite{WangWang2}.

Similar results and upper bounds to some of those above are known to hold in a related problem concerning proper total colourings
$c:E\cup V\to\{1,2,\ldots,k\}$ guaranteeing not only that $S(u)\neq S(v)$ for $uv\in E$, but also that
the sums of elements in $S(u)$ and $S(v)$ differ, see e.g.~\cite{PilsniakWozniak_total} (and~\cite{BonamyPrzybylo,FlandrinMPSW,Przybylo_CN_1,Przybylo_CN_2} for edge version of the same concept).
The least $k$ guaranteeing existence of such colouring
is denoted by $\chi''_{\sum}(G)$.
The definition of this, as well as of other graph invariants discussed above 
was motivated by a line of earlier research initiated by
Chartrand, Erd\H{o}s and Oellermann~\cite{ChartrandErdosOellermann}
and Chartrand et al.~\cite{Chartrand},
where a key seminal concept of the discipline, so called \emph{irregularity strength} of graphs was defined,
see e.g.~\cite{Aigner,KalKarPf,Nierhoff} for a few important results concerning this.
Though the property within the definition of $\chi''_{\sum}(G)$ is much stronger than the corresponding one required in the case of $\chi''_a(G)$,
and obviously $\chi''_a(G)\leq \chi''_{\sum}(G)$, no examples of graphs are known (to me) 
for which this inequality is not an equality in fact.
In this paper we wish to provide an evidence that indeed the difference between the two properties
is very significant. In order to bring out this we however need to extend the notions and definitions above
towards distinguishing not only neighbours, but also vertices at a greater, yet limited distance from each other.

For any fixed positive integer $r$,
vertices $u,v$ of $G$ shall be called \emph{$r$-neighbours} (or \emph{$r$-adjacent})
if $1\leq d(u,v)\leq r$, where $d(u,v)$ denotes the distance of $u$ and $v$ in $G$.
Similarly as in the concept of distant chromatic numbers (see~\cite{DistChrSurvey} for a survey of this topic),
the least number of colours in a proper total colouring $c$ of $G$ such that $S_c(u)\neq S_c(v)$
for every pair of 
$r$-neighbours $u,v\in V$,
so-called \emph{$r$-distant set distinguishing total colouring}
(or \emph{$r$-adjacent strong total colouring}),
shall be called the \emph{$r$-distant set distinguishing total number} or
\emph{$r$-adjacent strong total chromatic number}, and denoted by $\chi''_{a,r}(G)$.
This has already been considered (with different notation used) in several papers,
e.g. first in~\cite{ZhangDistant_total}, where mainly the case of $r=2$ was considered for several graph classes,
or under the name of the \emph{$d$-strong total chromatic number}
in~\cite{KemnitzMarangio_total}, where the parameter was settled for paths and cycles (for $r\leq 23$),
and a few 
bounds were given, in particular for cycles (in the remaining cases) and circulant graphs.
One of the main contribution of that paper relies however on disproving a general conjecture from~\cite{ZhangDistant_total}
concerning an upper bound for the value of the investigated graph invariant
via analysis of the cycles and circulant graphs, see also~\cite{Distsnt_total_Hatami_supplementary_1,ZhangDistant_total,Distsnt_total_Hatami_supplementary_3}
for some related supplementary results
(and cf.~\cite{Akbari,KemnitzMarangio,MocSotak,Przybylo_distant_Hatami_edge,Tian,ZhangDistant}
for results on edge correspondent of the same concept).
As for the known general upper bounds, it was proved in~\cite{Distsnt_total_Hatami_supplementary_general}
that if $r\geq 2$, then $\chi''_{a,r}(G)\leq 4\Delta^{r+2}$ if $3\leq \Delta\leq 4$, $\chi''_{a,r}(G)\leq 2\Delta^{r+2}$ if $5\leq \Delta\leq 6$
and $\chi''_{a,r}(G)\leq \Delta^{r+2}$ if $\Delta\geq 7$
for every graph $G$ with maximum degree $\Delta$.
We wish to strengthen these upper bounds significantly at a small, yet unavoidable cost.
Namely
we believe that
in order to distinguish not only neighbours, but also vertices at distance at most $r$ by colour pallets,
it still suffices to use $\Delta+C$ (instead of $\Omega(\Delta^{r+2})$) colours, where $C$ is some constant dependent on $r$,
if only the minimum degree of $G$ is greater than some constant,
linearly dependent on $r$, cf. the next section,
where we specify our main two results supporting 
our suspicion.
Note that if we wished to achieve the same using sums instead of sets,
one may quite easily construct a family of graphs with arbitrarily large minimum and maximum degrees 
which would require using (at least) $\Omega(\Delta^{r-1})$ colours,
cf. e.g.~\cite{Przybylo_distant} and the construction of graphs in Section~\ref{Construction_total_section}.
By our research we thus also want to expose
the leading impact
of the required \emph{properness} of the total colourings investigated
on the number of colours needed to distinguish
vertices by their incident sets. 
Indeed, we shall reveal that the number of these colours usually does not increase significantly
even if, as it might happen in our case, for every vertex $v$ of $G$,
the number of vertices from whose associated sets we need to distinguish
$S(v)$ 
grows considerably
(from at most $\Delta$ in the neighbourhood of $v$ possibly to roughly $\Delta^r$ in its $r$-neighbourhood).
The same could not hold in case of distinguishing by sums, see e.g.~\cite{Przybylo_distant}.
%
%
See also e.g. \cite{Hatami,Przybylo_distant_Hatami_edge,Zhang} for further motivation of our research.


\section{Results and Tools}

We pose the following conjecture.
\begin{conjecture}\label{Przybylo_main_conjecture_total}
For each positive integer $r$ there exist constants $\delta_0$ and $C$ such that
$$\chi''_{a,r}(G) \leq \Delta(G)+C$$
for every graph with $\delta(G) \geq \delta_0$.
\end{conjecture}
In fact in Section~\ref{Construction_total_section} we present a construction of a family of graphs
exemplifying that in order 
an upper bound for $\chi''_{a,r}(G)$ of the form
$\Delta+C$, or even $(1+o(1))\Delta$ to hold,
we cannot avoid some assumption on graphs considered,
e.g., that
the minimum
degree of $G$ is at least roughly (up to a small additive constant) $r$, see Observation~\ref{minimium_degree_obs_total}.
On the other hand, we believe that the conjecture above should hold with $\delta_0$ very close to $r$,
though we do not specify it explicitly.
The following two main results of this paper support this conjecture.
\begin{theorem}\label{first_attempt_total_theorem}
For every $r\geq 2$ there exists $\Delta_0$ such that for each graph $G$ of maximum degree $\Delta\geq\Delta_0$ with 
$\delta(G)\geq r+2$, 
$$\chi''_{a,r}(G)\leq \left\lceil\frac{\Delta}{\ln^4\Delta}\right\rceil (\lceil\ln^4\Delta\rceil+5\lceil\ln^3\Delta\rceil+2),$$
hence for all graphs $G$ of maximum degree $\Delta$ with $\delta(G)\geq r+2$,
$$\chi''_{a,r}(G) = (1+o(1))\Delta.$$
\end{theorem}
For every fixed $r$, we then strengthen the thesis, proving a desired upper bound of the form $\Delta+const.$
under additional assumption
that the minimum degree $\delta$ is not smaller than some function of $\Delta$,
however
disproportion between these two might be arbitrary large.
This e.g. encompasses the family of regular graphs (for which Conjecture~\ref{Przybylo_main_conjecture_total} thus holds).
\begin{theorem}\label{przybylo_main_total_bound}
For every positive $\varepsilon\leq 1$ and a positive integer $r$,
there exist $\Delta_0$ and a constant 
$C_2$ such that:
$$\chi''_{a,r}(G)\leq\Delta(G)+C_1+C_2$$
for every graph $G$ 
with $\delta(G)\geq \varepsilon\Delta(G)$ and $\Delta(G)\geq \Delta_0$,
where $C_1$ is the constant from Theorem~\ref{MolloyReedTh} and $C_2\leq \varepsilon^{-2}(7r+170)+r+5$.
\end{theorem}
The proof of this fact is inspired by the approaches in~\cite{CokerJohanson,Hatami,Przybylo_distant_Hatami_edge}.
%
%

We shall exploit a few tools of the probabilistic method, 
in particular
the Lov\'asz Local Lemma, see e.g.~\cite{AlonSpencer}, 
combined with the Chernoff Bound, see e.g.~\cite{JansonLuczakRucinski}  
(Th. 2.1, page 26)
and Talagrand's Inequality, see e.g.~\cite{MolloyReed}. 
\begin{theorem}[\textbf{The Local Lemma}]
\label{LLL-symmetric}
Let $A_1,A_2,\ldots,A_n$ be events in an arbitrary pro\-ba\-bi\-li\-ty space.
Suppose that each event $A_i$ is mutually independent of a set of all the other
events $A_j$ but at most $D$, and that ${\rm \emph{\textbf{Pr}}}(A_i)\leq p$ for all $1\leq i \leq n$. If
$$ ep(D+1) \leq 1,$$
then $ {\rm \emph{\textbf{Pr}}}\left(\bigcap_{i=1}^n\overline{A_i}\right)>0$.
\end{theorem}
\begin{theorem}[\textbf{Chernoff Bound}]\label{ChernofBoundTh}
For any $0\leq t\leq np$:
$${\rm\emph{\textbf{Pr}}}({\rm BIN}(n,p)>np+t)<e^{-\frac{t^2}{3np}}~~{and}~~{\rm\emph{\textbf{Pr}}}({\rm BIN}(n,p)<np-t)<e^{-\frac{t^2}{2np}}\leq e^{-\frac{t^2}{3np}}$$
where ${\rm BIN}(n,p)$ is the sum of $n$ independent Bernoulli variables, each equal to $1$ with probability $p$ and $0$ otherwise.
\end{theorem}
\begin{theorem}[\textbf{Talagrand's Inequality}]\label{TalagrandsInequalityTotal}
Let $X$ be a non-negative random variable, not identically $0$, which is determined by
$l$ independent trials $T_1,\ldots,T_l$, and satisfying the following for some $c,k>0$:
\begin{itemize}
\item[1.] changing the outcome of any one trial can affect $X$ by at most $c$, and
\item[2.] for any $s$, if $X\geq s$ then there is a set of at most $ks$ trials whose outcomes certify that $X\geq s$,
\end{itemize}
then for any $0\leq t\leq \mathbf{E}(X)$,
$${\rm\emph{\textbf{Pr}}}(|X-{\rm\emph{\textbf{E}}}(X)|>t+60c\sqrt{k{\rm\emph{\textbf{E}}}(X)})\leq 4e^{-\frac{t^2}{8c^2k\mathbf{E}(X)}}.$$
\end{theorem}

Note that e.g. knowing only an upper bound $\mathbf{E}(X)\leq h$ (instead of the exact value of $\mathbf{E}(X)$)
we may still use Talagrand's Inequality in order to upper-bound the probability that $X$ is large.
It is sufficient to apply Theorem~\ref{TalagrandsInequalityTotal} above
to the variable $Y=X+h-\mathbf{E}(X)$ with $\mathbf{E}(Y)=h$, to obtain the following
provided that the assumptions of Theorem~\ref{TalagrandsInequalityTotal}
hold for $X$ (and $t\leq h$):
$$\mathbf{Pr}(X>h+t+60c\sqrt{kh}) \leq \mathbf{Pr}(Y>h+t+60c\sqrt{kh}) \leq 4e^{-\frac{t^2}{8c^2kh}}.$$
Analogously, in the case of the Chernoff Bound, if $X$ is a sum of $n\leq k$ (where $k$ does not have to be an integer)
random independent Bernoulli variables, each equal to $1$ with probability $p\leq q$,
then $\mathbf{Pr}(X>kq+t)\leq e^{-\frac{t^2}{3kq}}$ (for $t\leq\lfloor k\rfloor q$),
and similarly, for $n\geq k\geq 2$ and $p\geq q$, $\mathbf{Pr}(X<kq-t)\leq e^{-\frac{t^2}{2\lceil k\rceil q}}\leq e^{-\frac{t^2}{3kq}}$ (for $t\leq\lceil k\rceil q$).
(It is sufficient to consider the variable $Y$ with binomial distribution ${\rm BIN}(k,\lfloor q\rfloor)$ or ${\rm BIN}(k,\lceil q\rceil)$,
respectively.)

We use the first two of these tools in the following section in order to prove Theorem~\ref{first_attempt_total_theorem}.
Talagrand's Inequality shall be used then in Section~\ref{SectionMainTotalProof} to obtain a strengthening of the thesis,
cf. Theorem~\ref{przybylo_main_total_bound}.
We shall also use the following well known inequalities a few times:
\begin{equation}\label{Choose_inequalities_total}
\left(\frac{a}{b}\right)^b\leq{a\choose b}
\leq\frac{a^b}{b!}\leq \left(\frac{ea}{b}\right)^b.
\end{equation}
Within the constructions below, some colours shall be removed from randomly selected edges.
In each such case we shall obtain so called \emph{partial colouring} $c$ of the considered graph, i.e.,
a mapping 
assigning colours only to some part of the edges and vertices of $G$.
For these, the colour pallets shall be understood (almost) the same as in~(\ref{totalScv}) (or~(\ref{edgeScv})),
but only the coloured edges and vertices shall be counted in the corresponding sets.

\section{Proof of Theorem~\ref{first_attempt_total_theorem}}
Let $r\geq 2$ be a fixed integer and $G=(V,E)$ be a graph with minimum degree $\delta(G)\geq r+2$.
Whenever needed we shall assume that its maximum degree $\Delta$ is sufficiently large.
Set $t=\lceil\frac{\Delta}{\ln^4\Delta}\rceil$ and $k=(\lceil\ln^4\Delta\rceil+5\lceil\ln^3\Delta\rceil+2)t$.
Partition the set of available colours $1,2,\ldots,k$ into $t$ (disjoint) subsets $C_1,C_2,\ldots,C_t$ of equal size.
Each $C_i$ we also partition into two subsets $C'_i$ and $C''_i$ with
$|C'_i|=\lceil\ln^4\Delta\rceil+\lceil\ln^3\Delta\rceil+1$
and $|C''_i|=4\lceil\ln^3\Delta\rceil+1$ for $i=1,\ldots,t$.

Let $c:V\to\{1,2,\ldots,k\}$ be any proper vertex colouring of $G$.
We shall extend it to a special total colouring of $G$,
which shall assure no conflicts between vertices of small degree
in all further steps of the construction. Then we shall randomly recolour a small portion
of the edges so that also $r$-neighbours of larger degrees become distinguished as well.

\subsection{Initial Total Colouring and vertices of small degrees}
For every edge $uv$ with $c(u)\in C_i$ and $c(v)\in C_j$ let us randomly and independently choose an integer
$\{1,2,\ldots,t\}\smallsetminus\{i,j\}$, each with equal probability. We denote this integer by $q(e)$ for every edge $e\in E$.
Note that the obtained auxiliary edge labelling $q$ does not have to be proper,
but we shall require in the further part of the construction that every $e\in E$ has always assigned a colour from $C_{q(e)}$
(by the choice of $q$ this shall in particular guarantee no colour conflicts between vertices and their incident edges).
Now for every $v\in V$, set $q(v)=i$ iff $c(v)\in C_i$ in order to make $q$ a \emph{total} labelling of $G$.
For every pair $u,v$ of $r$-neighbours with $d(u)=d(v)=d\leq \ln^3\Delta$ in $G$,
denote by
\begin{itemize}
\item $A_{\{u,v\}}$ - the event that $S_q(u)=S_q(v)$
\end{itemize}
(where the set $S_q(v)$ is understood as in~(\ref{totalScv})).
Then, since $u$ and $v$ might have at most one common incident edge (and $d\geq r+2$),
\begin{equation}\label{AuvProbability}
\mathbf{Pr}(A_{\{u,v\}})\leq \left(\frac{d+1}{t-2}\right)^{d-1} \leq \left(\frac{\ln^8\Delta}{\Delta}\right)^{r+1}.
\end{equation}
(as for every edge $e$ incident with $v$ we must have $q(e)\in S_q(u)$ if $A_{\{u,v\}}$ holds).
We shall also need a well distribution of of colours around every vertex.
Thus for every $v\in V$ and $i\in\{1,2,\ldots,t\}$ let us denote by
\begin{itemize}
\item $A_{v,i}$ - the event that more than $\ln^4\Delta+\ln^3\Delta$ edges $e$ incident with $v$ belong to $q^{-1}(i)$
(i.e., $q(e)=i$).
\end{itemize}
As $d(v)\leq \Delta$ and $\mathbf{Pr}(q(e)=i)\leq\frac{1}{t-2}\leq \frac{\ln^4\Delta+\frac{1}{2}\ln^3\Delta}{\Delta}$ for every $e\in E$,
then by the Chernoff Bound,
\begin{equation}\label{AviProbability}
\mathbf{Pr}(A_{v,i})\leq e^{-\frac{\frac{1}{4}\ln^6\Delta}{3(\ln^4\Delta+\frac{1}{2}\ln^3\Delta)}} \leq \left(\frac{\ln^8\Delta}{\Delta}\right)^{r+1}.
\end{equation}
Note that every event $A_{v,i}$ is mutually independent of all other events of the both types
defined above, except these indexed by some vertex which is at distance at most one from $v$,
i.e., except at most $(\Delta+1)t+\Delta\ln^3\Delta\cdot\Delta^{r-1} \leq 2\Delta^r\ln^3\Delta$ such events.
Analogously, every $A_{u,v}$ is mutually independent of all other events of the both types
except at most $2(\ln^3\Delta+1)t+2(\ln^3\Delta+1)\ln^3\Delta\cdot\Delta^{r-1} \leq 2\Delta^r\ln^3\Delta$.
Thus by (\ref{AuvProbability}) and (\ref{AviProbability}), the Local Lemma implies that (with positive probability)
$q$ may be chosen so that:
\begin{itemize}
\item[$(1^\circ)$] $S_q(u) \neq S_q(v)$ for every pair $u,v$ of $r$-neighbours with $d(u)=d(v)\leq \ln^3\Delta$ in $G$;
\item[$(2^\circ)$] at most $\ln^4\Delta+\ln^3\Delta$ edges incident with $v$ belong to $q^{-1}(i)$ for every $v\in V$ and $i\in\{1,2,\ldots,t\}$.
\end{itemize}
Then for every $i$, by $(2^\circ)$ above, the edges in $q^{-1}(i)$ induce a subgraph of $G$ with maximum degree
at most $\ln^4\Delta+\ln^3\Delta$. By Vizing's Theorem we may colour the edges of such subgraph properly with the colours from $C'_i$
subsequently for $i=1,2,\ldots,t$. We denote by $c(e)$ the colour assigned to every edge $e\in E$ in this manner.
These, together with the initial colours $c(v)$ chosen for vertices, complete the construction of our special proper total colouring $c$ of $G$.

\subsection{Recolouring and vertices of large degrees}
Now we randomly and independently uncolour the edges, each with probability $\frac{2}{\ln\Delta}$.
Denote the (partial) colouring obtained by $c'$, and let $U_{c'}$ be the set of uncoloured edges.
Denote also by $U_{c'}(v)$ the set of uncoloured edges incident with $v$.
We shall show that with positive probability the $r$-neighbours of sufficiently large degree are then distinguished.
In order to be able to complete our total colouring with relatively small number of colours at the end,
we however argue that uncoloured edges are also likely to be well distributed in general and within particular colour classes beforehand.

Consider any vertex $v\in V$ with $d(v)=d\geq \ln^3\Delta$.
Then $\mathbf{E}(|U_{c'}(v)|)=\frac{2d}{\ln\Delta}$,
and thus by the Chernoff Bound,
\begin{equation}\label{A_v_total_estimation1}
\mathbf{Pr}\left(\left||U_{c'}(v)|-\frac{2d}{\ln\Delta}\right|>\frac{d}{\ln\Delta}\right)<2e^{-\frac{d}{6\ln\Delta}}\leq 2e^{-\frac{\ln^2\Delta}{6}}\leq \frac{1}{\Delta^{r+3}}.
\end{equation}
Analogously, by $(2^\circ)$ above, $\mathbf{E}(|U_{c'}(v)\cap q^{-1}(i)|)\leq\frac{2(\ln^4\Delta+\ln^3\Delta)}{\ln\Delta} \leq 3\ln^3\Delta$,
and thus by the Chernoff Bound,
\begin{equation}\label{A_v_total_estimation2}
\mathbf{Pr}\left(|U_{c'}(v)\cap q^{-1}(i)|>4\ln^3\Delta\right)<e^{-\frac{(\ln^3\Delta)^2}{9\ln^3\Delta}}\leq \frac{1}{\Delta^{r+3}}.
\end{equation}
Denote by $B_{v,0}$ the event that $||U_{c'}(v)|-\frac{2d}{\ln\Delta}|>\frac{d}{\ln\Delta}$,
and for $i=1,\ldots,t$, denote by $B_{v,i}$ the event that $|U_{c'}(v)\cap q^{-1}(i)|>4\ln^3\Delta$ (for $d\geq \ln^3\Delta$).
For any $r$-neighbours $u,v$ which are of the same degree $d$ with $\ln^3\Delta\leq d\leq\Delta$ in $G$,
denote by $B_{\{u,v\}}$ the event that $|U_{c'}(u)|,|U_{c'}(v)|\in \left[\frac{d}{\ln\Delta},\frac{3d}{\ln\Delta}\right]$
and $S_{c'}(u)=S_{c'}(v)$.
Since $u$ and $v$ have at most one common incident edge, then:
\begin{eqnarray}
\mathbf{Pr}(B_{\{u,v\}}) &\leq& \mathbf{Pr}\left(S_{c'}(u)=S_{c'}(v)\wedge |U_{c'}(v)|\in \left[\frac{d}{\ln\Delta},\frac{3d}{\ln\Delta}\right]\right)\nonumber\\
&\leq& \mathbf{Pr}\left(S_{c'}(u)=S_{c'}(v)\left||U_{c'}(v)|\in \left[\frac{d}{\ln\Delta},\frac{3d}{\ln\Delta}\right]\right.\right)\nonumber\\
&\leq&
\sum_{j =\lceil\frac{d}{\ln\Delta}\rceil}^{\lfloor\frac{3d}{\ln\Delta}\rfloor}
\mathbf{Pr}\left(S_{c'}(u)=S_{c'}(v)\left||U_{c'}(v)|=j\right.\right)
\mathbf{Pr}\left(|U_{c'}(v)|=j\left||U_{c'}(v)|\in \left[\frac{d}{\ln\Delta},\frac{3d}{\ln\Delta}\right]\right.\right)\nonumber\\
&\leq&
\sum_{j =\lceil\frac{d}{\ln\Delta}\rceil}^{\lfloor\frac{3d}{\ln\Delta}\rfloor}
\left(\frac{2}{\ln\Delta}\right)^{j-1}
\mathbf{Pr}\left(|U_{c'}(v)|=j\left||U_{c'}(v)|\in \left[\frac{d}{\ln\Delta},\frac{3d}{\ln\Delta}\right]\right.\right)\nonumber\\
&\leq&
\left(\frac{2}{\ln\Delta}\right)^{\frac{d}{\ln\Delta}-1}
\sum_{j =\lceil\frac{d}{\ln\Delta}\rceil}^{\lfloor\frac{3d}{\ln\Delta}\rfloor}
\mathbf{Pr}\left(|U_{c'}(v)|=j\left||U_{c'}(v)|\in \left[\frac{d}{\ln\Delta},\frac{3d}{\ln\Delta}\right]\right.\right)\nonumber\\
&\leq& \left(\frac{2}{\ln\Delta}\right)^{\ln^2\Delta-1}\cdot 1
\leq \left(\frac{1}{e}\right)^{(r+3)\ln\Delta}
=\frac{1}{\Delta^{r+3}}.\label{B_uv_total_estimation}
\end{eqnarray}
Analogously as above, every event
$B_{v,i}$ 
or $B_{\{u,v\}}$ is mutually independent of all other events
of these forms indexed by vertices each of which is at distance at least
$2$ from both $v$ and $u$ (in the case of $B_{\{u,v\}}$), i.e., of all but at most
$2(t+2)(\Delta+1) \Delta^r \leq \Delta^{r+2}$ other events.
At the same time, by~(\ref{A_v_total_estimation1}),~(\ref{A_v_total_estimation2}) and~(\ref{B_uv_total_estimation}) each of these events occurs with probability at most $\frac{1}{\Delta^{r+3}}$.
By the 
Local Lemma, with positive probability the uncoloured edges could be chosen
so that none of the events $B_{v,i}$ and $B_{\{u,v\}}$ holds.
As the events of the forms $B_{v,0}$ or $B_{\{u,v\}}$ do not appear,
every pair of $r$-neighbours of the same degree $d\geq \ln^3\Delta$ is set distinguished in $G$ then.
Moreover, for every $i=1,2,\ldots,t$, by nonappearance of the events $B_{v,i}$,
the subgraph induced by the edges in $U_{c'}\cap q^{-1}(i)$ has maximum degree at most
$4\ln^3\Delta$, and thus can be coloured properly with
(so far spared) colours from $C''_i$.
Since we have used new colours to paint the uncoloured edges, the obtained total colouring $c''$ of $G$ is then proper,
as the new colours still belonged to $C_{q(e)}$ for every edge $e\in E$ (cf. the definition of $q$).
Using new colours also guarantees preservation of distinction within pairs of $r$-neighbours of large degrees.
On the other hand, if $u$ and $v$ are $r$-neighbours with $d(u)=d(v)\leq\ln^3\Delta$,
then $S_{c''}(u)\neq S_{c''}(v)$ by the condition $(1^\circ)$. Therefore, the total colouring $c''$ is $r$-distant set distinguishing.
The proof of Theorem~\ref{first_attempt_total_theorem} is thus completed.
$\blacksquare$


\section{Proof of Theorem~\ref{przybylo_main_total_bound}\label{SectionMainTotalProof}}

Fix any $\varepsilon\in (0,1]$ and a positive integer $r$.
Let $G=(V,E)$ be a graph with $\delta\geq \varepsilon\Delta$, where $\Delta$ and $\delta$
are its maximum and minimum degrees, respectively.
Whenever needed it shall be assumed that $\Delta$ is sufficiently large, i.e., we explicitly do not
specify $\Delta_0$.
We shall assign to edges and vertices the colours $1,2,\ldots,\Delta+C_1+C_2$, where
$C_1$ is the constant from Theorem~\ref{MolloyReedTh} and $C_2=\lfloor\varepsilon^{-2}(7r+170)+r+5\rfloor$, in two stages in order to obtain
an $r$-distant set distinguishing total colouring of $G$.

We shall start from a given total colouring of $G$, and in the first stage of our randomized construction,
similarly as in the previous proof, we shall uncolour some of the edges.
This time, we shall however admit only a constant number of uncoloured edges from every $E(v)$, $v\in V$.
Analogously as above, $r$-neighbours with sufficiently many uncoloured incident edges shall be very likely
to be set distinguished, but we shall not 
avoid appearances of
troublesome vertices with very few uncoloured incident colours.
Enhancing our argument with Talagrand's Inequality, we shall however be able to show that
these should be rare and well distributed.
In the second stage we shall additionally uncolour a few edges incident with every troublesome vertex,
thus distinguishing it form its $r$-neighbours.
In order not to spoil distinctions from the first stage, these shall be required to be strong enough,
i.e., we shall require the sets of (not troublesome) vertices to differ in sufficiently many elements
from the sets associated with their $r$-neighbours of the same degree,
see \emph{(b)} in Claim~\ref{FinalClaimSt1_total} on the symmetric difference of such sets below.
A constant number of extra colour shall be used to complete the total colouring at the end.

\subsection{Stage One}
Fix any proper total colouring $c_0:V\cup E\to\{1,2,\ldots,\Delta+C_1\}$ of $G$, guaranteed by Theorem~\ref{MolloyReedTh}.
From this base colouring we create a new (partial) total colouring $c$ of $G$ in two steps as follows:
\begin{itemize}
\item uncolour each edge $e\in E$ independently with probability $\frac{5r+90}{\varepsilon^2\Delta}$;
denote the set of uncoloured edges in this step by $U$, and the subset of these incident with any given $v\in V$ by $U(v)$;
\item and then, for every vertex with more than $\varepsilon^{-2}(7r+170)$ incident edges
uncoloured in the step above, we recover the removed colours of all its incident edges;
call the corresponding vertices \emph{recovered}.
\end{itemize}
We shall also denote by $U_c(v)$ the set of edges incident with $v\in V$ which are not coloured under $c$
(where $c$ is the resulting colouring after the two steps described above).
Let $L$ be the set of all (\emph{troublesome}) vertices $v\in V$ with $|U_c(v)|< 3r+15$.
Note that
\begin{equation}\label{Lincludedinsum}
L\subseteq R\cup L_U\cup L_R,
\end{equation}
where:

\begin{itemize}
\item $R$ is the set of all recovered vertices;
\item $L_U$ is the set of vertices $v$ with $|U(v)|<3r+15$; 
\item $L_R$ is the set of all vertices $v$ adjacent with some vertex $u\in R$ such that $uv\in U$.
\end{itemize}

We first provide a small, but useful technical observation, repeatedly exploited in the further part of the argument.
\begin{observation}\label{two_technical_inequalities_total}
For every $\varepsilon \in (0,1]$ and a positive integer $r$,
\begin{eqnarray}
e^{10(1-\varepsilon^{-1})} &\leq& e^{10(1-\varepsilon^{-1})}\varepsilon^{-2} \leq \varepsilon^2,\label{EpsilonIneq}\\
e^{-\frac{4r+150}{15}} &\leq& e^{-\frac{4r+150}{15}}(5r+90) \leq \frac{1}{100r}.\label{rIneq}
\end{eqnarray}
\end{observation}
\begin{pf}
The left-hand side inequalities above are straightforward. We shall justify the remaining ones.
The one from (\ref{EpsilonIneq}) is equivalent to the fact that
$$f(\varepsilon)= 10(1-\varepsilon^{-1})-4\ln\varepsilon\leq 0,$$
which is true, as
$$f'(\varepsilon)=\frac{4(2.5-\varepsilon)}{\varepsilon^2}$$
(hence $f$ is increasing for $\varepsilon\in(0,1]$) and $f(1)=0$.

The one from (\ref{rIneq}) is in turn equivalent to the inequality:
\begin{equation}\label{EquivGIneq}
\frac{4r+150}{15} - \ln[100r(5r+90)]\geq 0,
\end{equation}
but since $100r(5r+90) 
= (30r+90)^2 -(20r-90)^2 \leq (30r+90)^2$ (for $r\geq 1$),
to prove (\ref{EquivGIneq}) it is then sufficient to observe that
$$g(r)=\frac{4r+150}{15} - \ln(30r+90)^2\geq 0$$
for $r\geq 1$, as
$$g'(r)=\frac{4}{15}-\frac{2}{r+3}
=\frac{4(r-4.5)}{15(r+3)}$$
(i.e., $g$ is decreasing for $r\in[1,4.5]$ and
increasing for $r\geq 4.5$) and
$g(4.5)\approx 0.37 > 0$.
$\blacksquare$
\end{pf}

\begin{claim}\label{main_total_claim}
For every vertex $v\in V$,
\begin{equation}\label{main_total_probability}
\mathbf{Pr}\left(|N(v)\cap L| > \frac{\varepsilon^2\Delta}{10r}\right) \leq \frac{1}{10\Delta^{r+4}}.
\end{equation}
\end{claim}

\begin{pf}
By (\ref{Lincludedinsum}) it is sufficient to prove the following three inequalities:
\begin{eqnarray}
\mathbf{Pr}\left(|N(v)\cap R| > \frac{\varepsilon^2\Delta}{30r}\right) &\leq& \frac{1}{30\Delta^{r+4}}\label{main_total_probability_1},\\
\mathbf{Pr}\left(|N(v)\cap L_U| > \frac{\varepsilon^2\Delta}{30r}\right) &\leq& \frac{1}{30\Delta^{r+4}}\label{main_total_probability_2},\\
\mathbf{Pr}\left(|N(v)\cap L_R| > \frac{\varepsilon^2\Delta}{30r}\right) &\leq& \frac{1}{30\Delta^{r+4}}\label{main_total_probability_3}.
\end{eqnarray}

Since for any given vertex $u\in V$, the random variable $|U(u)|$ has binomial distribution
with parameters $d(u)\leq\Delta$ and $\frac{5r+90}{\varepsilon^2\Delta}$,
by the Chernoff Bound,
\begin{eqnarray}
{\rm\textbf{Pr}}(u\in R) &=& {\rm\textbf{Pr}}(|U(u)|>\varepsilon^{-2}(7r+170))
\leq e^{-\frac{[\varepsilon^{-2}(7r+170)- \varepsilon^{-2}(5r+90)]^2}{3 \varepsilon^{-2}(5r+90)}}\nonumber\\
&=& e^{-\frac{\varepsilon^{-2}(2r+80)^2}{3(5r+90)}}
\leq e^{-\frac{2\varepsilon^{-2}(2r+80)}{15}}
= e^{-\frac{4r+160}{15}(\varepsilon^{-2}-1)}\cdot e^{-\frac{4r+160}{15}}\nonumber\\
&\leq& e^{10(1-\varepsilon^{-1})}\cdot e^{-\frac{4r+150}{15}}
\leq \varepsilon^2\cdot\frac{1}{100r},\label{vinRineq_total}
\end{eqnarray}
where the last inequality above holds by Observation~\ref{two_technical_inequalities_total}.
Analogously, as $d(u)\geq \varepsilon\Delta$, then by the Chernoff Bound,
\begin{eqnarray}
{\rm\textbf{Pr}}(u\in L_U) &=& {\rm\textbf{Pr}}(|U(u)|<3r+15) \leq {\rm\textbf{Pr}}(|U(u)|<\varepsilon^{-1}(3r+15))\nonumber\\
&\leq& e^{-\frac{[\varepsilon^{-1}(5r+90)- \varepsilon^{-1}(3r+15)]^2}{3 \varepsilon^{-1}(5r+90)}}
= e^{-\frac{\varepsilon^{-1}(2r+75)^2}{3(5r+90)}}
\leq e^{-\frac{2\varepsilon^{-1}(2r+75)}{15}}\nonumber\\
&=& e^{-\frac{4r+150}{15}(\varepsilon^{-1}-1)}\cdot e^{-\frac{4r+150}{15}}
\leq e^{10(1-\varepsilon^{-1})}\cdot e^{-\frac{4r+150}{15}}
\leq \varepsilon^2\cdot\frac{1}{100r},\label{vinL_Uineq_total}
\end{eqnarray}
again due to Observation~\ref{two_technical_inequalities_total}.
Additionally, if $w$ is any neighbour of $u$, then by the Chernoff Bound (similarly as in (\ref{vinRineq_total})),
as $|E(u)\smallsetminus\{uw\}|\leq \Delta-1\leq \Delta$,
\begin{eqnarray}
{\rm\textbf{Pr}}(w\in R|uw\in U)
&\leq& e^{-\frac{[\varepsilon^{-2}(7r+170)-1- \varepsilon^{-2}(5r+90)]^2}{3 \varepsilon^{-2}(5r+90)}}
\leq e^{-\frac{\varepsilon^{-2}(2r+79)^2}{3(5r+90)}}, 
\nonumber
\end{eqnarray}
hence,
\begin{eqnarray}
{\rm\textbf{Pr}}(w\in R\wedge uw\in U)&=& {\rm\textbf{Pr}}(w\in R|uw\in U)\cdot {\rm\textbf{Pr}}(uw\in U)
\leq e^{-\frac{\varepsilon^{-2}(2r+79)^2}{3(5r+90)}} \cdot \frac{\varepsilon^{-2}(5r+90)}{\Delta}.\nonumber
\end{eqnarray}
Consequently,
\begin{eqnarray}
{\rm\textbf{Pr}}(u\in L_R)&\leq& \Delta\cdot e^{-\frac{\varepsilon^{-2}(2r+79)^2}{3(5r+90)}} \cdot \frac{\varepsilon^{-2}(5r+90)}{\Delta}
\leq e^{-\frac{2\varepsilon^{-2}(2r+79)}{15}} \cdot \varepsilon^{-2}(5r+90)\nonumber\\
&=& e^{-\frac{4r+158}{15}(\varepsilon^{-2}-1)} \varepsilon^{-2} \cdot e^{-\frac{4r+158}{15}} (5r+90)\nonumber\\
&\leq& e^{10(1-\varepsilon^{-1})}\varepsilon^{-2}\cdot e^{-\frac{4r+150}{15}}(5r+90)
\leq \varepsilon^2\cdot\frac{1}{100r},
\label{vinL_Rineq_total}
\end{eqnarray}
where the last inequality follows by Observation~\ref{two_technical_inequalities_total}.

With inequalities (\ref{vinRineq_total})-(\ref{vinL_Rineq_total}) at hand we shall now use Talagrand's Inequality
in order to prove (\ref{main_total_probability_1})-(\ref{main_total_probability_3}).
For every vertex $v\in V$, by (\ref{vinRineq_total}), ${\rm\textbf{E}}(|N(v)\cap R|)\leq \frac{\varepsilon^2\Delta}{100r}$.
Note also that the random variable $|N(v)\cap R|$ 
is being determined by the outcomes
of Bernoulli trials associated with all edges incident with neighbours of $v$,
each of which sets down whether a colour is removed from
a given edge in the first step of the construction or not (i.e., whether the edge belongs to $U$ or not).
Moreover, changing the outcome of any one such trial may affect $|N(v)\cap R|$ 
by at most $2$,
while the fact that $|N(v)\cap R|\geq s$ 
can be certified by the outcomes of at most
$(\varepsilon^{-2}(7r+170)+1)s$ trials. Therefore,
\begin{eqnarray*}
{\rm\textbf{Pr}}\left(|N(v)\cap R| > \frac{\varepsilon^2\Delta}{30r} \right)
&\leq& {\rm\textbf{Pr}}\left(|N(v)\cap R| > \frac{\varepsilon^2\Delta}{100r}+\frac{\varepsilon^2\Delta}{100r}+120\sqrt{(\varepsilon^{-2}(7r+170)+1)\frac{\varepsilon^2\Delta}{100r}}\right)\\
&\leq& 4e^{-\frac{(\frac{\varepsilon^2\Delta}{100r})^2}{8\cdot 4\cdot (\varepsilon^{-2}(7r+170)+1)\frac{\varepsilon^2\Delta}{100r}}}< \frac{1}{30\Delta^{r+4}}
\end{eqnarray*}
(for $\Delta$ sufficiently large), hence inequality (\ref{main_total_probability_1}) holds.

Analogously, for each vertex $v$,
by (\ref{vinL_Rineq_total}), ${\rm\textbf{E}}(|N(v)\cap L_R|)\leq \frac{\varepsilon^2\Delta}{100r}$.
This time the variable $|N(v)\cap L_R|$ is determined by the outcomes
of Bernoulli trials associated with all edges incident with neighbours of $v$ or their neighbours
(determining belongingness of these in $U$).
Changing the outcome of any such trial may affect $|N(v)\cap L_R|$ by at most $2(\varepsilon^{-2}(7r+170)+1)$,
while the fact that $|N(v)\cap L_R|\geq s$ can be certified by the outcomes of at most
$(\varepsilon^{-2}(7r+170)+1)s$ trials. Thus,
\begin{eqnarray*}
&&{\rm\textbf{Pr}}\left(|N(v)\cap L_R| > \frac{\varepsilon^2\Delta}{30r} \right)\\
&\leq& {\rm\textbf{Pr}}\left(|N(v)\cap L_R| > \frac{\varepsilon^2\Delta}{100r}+\frac{\varepsilon^2\Delta}{100r} + 120(\varepsilon^{-2}(7r+170)+1) \sqrt{(\varepsilon^{-2}(7r+170)+1)\frac{\varepsilon^2\Delta}{100r}}\right)\\
&\leq& 4e^{-\frac{(\frac{\varepsilon^2\Delta}{100r})^2}{8\cdot 4\cdot (\varepsilon^{-2}(7r+170)+1)^3 \frac{\varepsilon^2\Delta}{100r}}}
< \frac{1}{30\Delta^{r+4}},
\end{eqnarray*}
i.e., inequality (\ref{main_total_probability_3}) holds.

We shall be a slightly more careful with inequality (\ref{main_total_probability_2})
though.
In fact we shall bound the probability that the random variable $X=d(v)-|N(v)\cap L_U|$ is small instead.
By~(\ref{vinL_Uineq_total}) we have:
$$\Delta \geq {\rm\textbf{E}}(X) = d(v)-{\rm\textbf{E}}(|N(v)\cap L_U|) \geq d(v)-\frac{\varepsilon^2\Delta}{100r}.$$
As previously, the variable $X$
is determined by the outcomes
of Bernoulli trials associated with all edges incident with neighbours of $v$,
changing the outcome of each of which may affect $X$ by at most $2$.
As the fact that $X\geq s$ can be certified by the outcomes of at most
$(3r+15)s$ trials, by Talagrand's Inequality we obtain:
\begin{eqnarray*}
&&{\rm\textbf{Pr}}\left(d(v)-|N(v)\cap L_U|<d(v) - \frac{\varepsilon^2\Delta}{30r}\right)\\
&\leq&{\rm\textbf{Pr}}\left(X<\left(d(v)-\frac{\varepsilon^2\Delta}{100r}\right)-\frac{\varepsilon^2\Delta}{100r} - 120\sqrt{(3r+15)\Delta}\right)\\
&\leq&{\rm\textbf{Pr}}\left(|X-{\rm\textbf{E}}(X)|>\frac{\varepsilon^2\Delta}{100r} + 120\sqrt{(3r+15){\rm\textbf{E}}(X)}\right)\\
&\leq& 4e^{-\frac{(\frac{\varepsilon^2\Delta}{100r})^2}{8\cdot 2^2(3r+15)\mathbf{E}(X)}}
\leq 4e^{-\frac{(\frac{\varepsilon^2\Delta}{100r})^2}{8\cdot 2^2(3r+15)\Delta}}
\leq \frac{1}{30\Delta^{r+4}},
\end{eqnarray*}
(for $\Delta$ sufficiently large),
hence (\ref{main_total_probability_2}) follows.
$\blacksquare$
\end{pf}

For any sets $A$ and $B$, let $A\triangle B:=(A\smallsetminus B)\cup(B\smallsetminus A)$ denote their \emph{symmetric difference}.

\begin{claim}\label{Claim_symm_diff_total}
For every pair of $r$-neighbours $u,v$ with $d(u)=d(v)$ and $1\leq d(u,v)\leq r$ in $G$,
$${\rm\emph{\textbf{Pr}}}(u\notin L \wedge |S_c(u)\triangle S_c(v)|<2r+10)<\frac{1}{10\Delta^{r+4}}.$$
\end{claim}

\begin{pf}
Let $u,v\in V$, $d(u)=d(v)$ and $1\leq d(u,v)\leq r$.
Note that $u\notin L$, i.e. $|U_c(u)|\geq 3r+15$ implies in particular that $3r+15\leq |U(u)|\leq \varepsilon^{-2}(7r+170)$, hence
\begin{eqnarray}
& & {\rm\textbf{Pr}}(u\notin L \wedge |S_c(u)\triangle S_c(v)|<2r+10) \nonumber\\
&=& {\rm\textbf{Pr}}(|U_c(u)|\geq 3r+15 \wedge 3r+15\leq |U(u)|\leq \varepsilon^{-2}(7r+170) \wedge |S_c(u)\triangle S_c(v)|\leq 2r+9) \nonumber\\
&\leq& {\rm\textbf{Pr}}(|U_c(u)|\geq 3r+15 \wedge |S_c(u)\triangle S_c(v)|\leq 2r+9~|~3r+15\leq |U(u)|\leq \varepsilon^{-2}(7r+170)). \nonumber
\end{eqnarray}
Now, since the random uncolourings of the edges within our process are independent,
suppose we first perform the corresponding experiments (determining whether a given edge is uncoloured or not)
for the edges incident with $u$, and that $3r+15\leq |U(u)|\leq \varepsilon^{-2}(7r+170)$.
After executing only these, there are at least $|U(u)|$ (or $|U(u)|-1$ if $uv\in E$ and $uv$ was uncoloured) elements incident with $v$
(i.e., edges incident with $v$ or $v$ itself)
whose colours do not belong to the pallet of $u$.
Out of these choose $|U(u)|$ ($|U(u)|-1$, resp.) coloured with the least integers
and denote them by $E'_v$ (where this set might include $v$).
Since we wish to have $|U_c(u)|\geq 3r+15$ within the investigated conditional event,
at most $|U(u)|-3r-15$ edges in $U(u)$ might have their
colours recovered 
eventually (some or all of which might be assigned to the edges or vertex in $E'_v$).
As $|E'_v|-(|U(u)|-3r-15)\geq 3r+14$,
in order to have 
$|S_c(u)\triangle S_c(v)|\leq 2r+9$ at the end,
still at least $r+5$ edges from $E'_v$ must be uncoloured in our random process.
Since 
$|U(u)| \leq \varepsilon^{-2}(7r+170)$, hence $|E'_v| \leq \varepsilon^{-2}(7r+170)$, we finally obtain that:
\begin{eqnarray}
{\rm\textbf{Pr}}(u\notin L \wedge |S_c(u)\triangle S_c(v)|<2r+10) 
&\leq& {\lfloor\varepsilon^{-2}(7r+170)\rfloor \choose r+5} \left(\frac{\varepsilon^{-2}(5r+90)}{\Delta}\right)^{r+5} \nonumber\\
&<& \frac{1}{10\Delta^{r+4}} \nonumber
\end{eqnarray}
(for $\Delta$ sufficiently large).
$\blacksquare$
\end{pf}

\begin{claim}\label{FinalClaimSt1_total}
With positive probability, we have:
\begin{itemize}
\item[(a)] $|N(v)\cap L| \leq \frac{\varepsilon^2\Delta}{10r}$
for every vertex $v\in V$, and
\item[(b)] $|S_c(u)\triangle S_c(v)|\geq 2r+10$ for every pair of $r$-neighbours $u,v$
with $d(u)=d(v)$ and $u\notin L$ (or $v\notin L$).
\end{itemize}
\end{claim}

\begin{pf}
For any vertex $v\in V$, let $A_v$ be the event that $|N(v)\cap L| > \frac{\varepsilon^2\Delta}{10r}$.
For every $u,v\in V$ with $d(u)=d(v)$ and $1\leq d(u,v)\leq r$, let $B_{u,v}$ be the event
that $u\notin L$ and $|S_c(u)\triangle S_c(v)|<2r+10$
(note that we must distinguish $B_{u,v}$ from $B_{v,u}$ according to this definition).
In order to apply Theorem~\ref{LLL-symmetric}, let us first
notice that every event $A_v$ is mutually independent of all events $A_{v'}$ and $B_{u,w}$
with $d(v,v')>5$, $d(v,u)>4$ and $d(v,w)>4$, i.e., of all other events of these forms
but at most $\Delta^5+\Delta^4\cdot\Delta^r\cdot 2 \leq 3\Delta^{r+4}$.
Similarly, each event $B_{u,v}$ is mutually independent of all events $A_w$ and $B_{u',v'}$
with $d(u,w)>4$, $d(v,w)>4$ and $d(u,u')>3$, $d(u,v')>3$, $d(v,u')>3$, $d(v,v')>3$,
i.e., of all other events of these forms
but at most $2\Delta^4+2\cdot\Delta^3\cdot\Delta^r\cdot 2 \leq 3\Delta^{r+4}$.
As by Claims~\ref{main_total_claim} and~\ref{Claim_symm_diff_total},
each of these events occurs with probability at most $\frac{1}{10\Delta^{r+4}}$,
the thesis follows by the Lov\'asz Local Lemma.
$\blacksquare$
\end{pf}

\subsection{Stage Two}

Suppose $c$ is a partial (proper) total colouring constructed in Stage One and satisfying the thesis of Claim~\ref{FinalClaimSt1_total} above.
Prior completing this we still need to distinguish the troublesome vertices in $L$
from their $r$-neighbours of the same degrees.
For this aim we now independently for every such vertex $v\in L$
randomly uncolour its $r+5$
incident edges which were coloured under $c$ and joined $v$ with vertices outside $L$.
The partial colouring obtained is then denoted by $c'$, and the set of uncoloured within this stage edges by $U'$.
We also denote the set of edges from $U'$ which are incident with any given $v\in V$ by $U'(v)$.
In the following claim we shall argue that the colouring $c'$
does not have to influence the colour pallets of vertices outside $L$ significantly
(and thus the condition
\emph{(b)} of Claim~\ref{FinalClaimSt1_total} above shall suffice to keep them distinguished from their $r$-neighbours of the same degrees),
and at the same time $c'$ might guarantee distinction between vertices in $L$.

\begin{claim}\label{mainSecondStageClaim_total}
With positive probability,
\begin{itemize}
\item[(i)] $|U'(u)|\leq r+4$ for every vertex $u\in V\smallsetminus L$, and
\item[(ii)] $S_{c'}(u)\neq S_{c'}(v)$ for every $u,v\in L$ with $d(u)=d(v)$ and $1\leq d(u,v)\leq r$.
\end{itemize}
\end{claim}
\begin{pf}
Let us define two kinds of (bad) events corresponding to \emph{(i)} and \emph{(ii)}, respectively:
\begin{itemize}
\item for any $r+5$ neighbours $v_1,\ldots,v_{r+5}\in L$
  of a given vertex $u\in V\smallsetminus L$
  such that
  $uv_i\notin U_c(u)$ for $i=1,\ldots,r+5$,
let $A_{u,\{v_1,\ldots,v_{r+5}\}}$ be the event that $uv_i\in U'$ for $i=1,\ldots,r+5$;
\item for every vertices $u,v\in L$ with $d(u)=d(v)$ and $1\leq d(u,v)\leq r$, let $B_{\{u,v\}}$ be the event that $S_{c'}(u)=S_{c'}(v)$.
\end{itemize}

Since by Stage One, every vertex $v\in L$ has at least
$\varepsilon\Delta-\frac{\varepsilon^{2}\Delta}{10r}-(3r+15)
\geq \frac{4}{5}\varepsilon\Delta$
neighbours $w\notin L$
such that
$vw\notin U_c(v)$, then: 
\begin{eqnarray}
{\rm\textbf{Pr}}(A_{u,\{v_1,\ldots,v_{r+5}\}}) &\leq& \left(\frac{r+5}{\frac{4}{5}\varepsilon\Delta}\right)^{r+5} ~~ {\rm and} \label{Auvi_ineq_total}\\
{\rm\textbf{Pr}}(B_{\{u,v\}}) &\leq& \frac{1}{{\lceil\frac{4}{5}\varepsilon\Delta\rceil \choose r+5}}
\leq \frac{1}{\left(\frac{\frac{4}{5}\varepsilon\Delta}{r+5}\right)^{r+5}}
= \left(\frac{r+5}{\frac{4}{5}\varepsilon\Delta}\right)^{r+5}. \label{Buv_ineq_total}
\end{eqnarray}

Note also that every event $A_{u,\{v_1,\ldots,v_{r+5}\}}$ is mutually independent of all such events (of the both types)
except for these $A_{u',\{v'_1,\ldots,v'_{r+5}\}}$ and $B_{\{u',v'\}}$ for which we have
$\{v_1,\ldots,v_{r+5}\}\cap\{v'_1,\ldots,v'_{r+5}\}\neq \emptyset$ or $\{v_1,\ldots,v_{r+5}\}\cap\{u',v'\}\neq \emptyset$, resp.,
i.e., of all other events except at most $(r+5)\Delta {\lfloor\frac{\varepsilon^2\Delta}{10r}\rfloor \choose r+4} + (r+5)\Delta^{r}$.
Analogously, every event $B_{\{u,v\}}$ is mutually independent of all such events
except for these $A_{u',\{v_1,\ldots,v_{r+5}\}}$ and $B_{\{u',v'\}}$ for which we have
$\{u,v\}\cap\{v_1,\ldots,v_{r+5}\}\neq \emptyset$ or $\{u,v\}\cap\{u',v'\}\neq\emptyset$, resp.,
i.e., of all other events except at most $2\Delta {\lfloor\frac{\varepsilon^2\Delta}{10r}\rfloor \choose r+4} + 2\Delta^{r}
\leq (r+5)(\Delta {\lfloor\frac{\varepsilon^2\Delta}{10r}\rfloor \choose r+4} + \Delta^{r})$.
As
\begin{eqnarray}
e \left(\frac{r+5}{\frac{4}{5}\varepsilon\Delta}\right)^{r+5} [(r+5)(\Delta {\lfloor\frac{\varepsilon^2\Delta}{10r}\rfloor \choose r+4} + \Delta^{r})+1]
&\leq& 4\left(\frac{r+5}{\frac{4}{5}\varepsilon\Delta}\right)^{r+5}
(r+5)\Delta\frac{\left(\frac{\varepsilon^2\Delta}{10r}\right)^{r+4}}{(r+4)!} \nonumber\\
&\leq& \frac{5(r+5)^{r+6}}{(r+4)!(8r)^{r+4}}
< 1 \nonumber
\end{eqnarray}
(where the last inequality can be checked directly for $r=1$, while for $r\geq 2$ we obviously have:
$(r+5)^{r+4}<(8r)^{r+4}$ and $5(r+5)^2<(r+4)!$),
this together with inequalities (\ref{Auvi_ineq_total}) and (\ref{Buv_ineq_total})
suffices to apply the Local Lemma, and thus obtain the thesis.
$\blacksquare$
\end{pf}

Let $c'$ be a partial (proper) colouring of $G$ guaranteed to exist by Claim~\ref{mainSecondStageClaim_total}.
By \emph{(b)} from Claim~\ref{FinalClaimSt1_total} and \emph{(i)},\emph{(ii)} from Claim~\ref{mainSecondStageClaim_total},
all $r$-neighbours of the same degrees are distinguished under $c'$.
By our construction, every vertex $v\notin L$ is incident with at most
$\varepsilon^{-2}(7r+170) + (r+4)$ uncoloured edges under $c'$,
while every $v\in L$ is incident with no more than $(3r+15) + (r+5)\leq \varepsilon^{-2}(7r+170) + (r+4)$
such edges. By Vizing's Theorem, we thus need at most $\lfloor\varepsilon^{-2}(7r+170) + (r+4)\rfloor+1$
extra colours in order to extend $c'$ to a proper total colouring of $G$
(all together using at most $\Delta+C_1+\varepsilon^{-2}(7r+170) + r+5$ colours).
As these new colours shall not spoil distinction between the pallets associated to $r$-neighbours,
we have obtained a desired total colouring of $G$. The proof of Theorem~\ref{przybylo_main_total_bound} is thus completed.
$\blacksquare$

\section{Construction\label{Construction_total_section}}

We shall make use of so called \emph{undirected de Bruijn graph} of type $(t,k)$, denoted by $D_{t,k}$ and defined as follows.
The vertex set of $D_{t,k}$ is formed of all sequences of length $k$
the entries of which are taken from a fixed alphabet consisting of $t$ distinct letters.
The edges of such graph are in turn formed by any two distinct vertices $(a_1,\ldots,a_k)$ and $(b_1,\ldots,b_k)$ for which
either $a_i = b_{i+1}$ for $1\leq i \leq k - 1$, or $a_{i+1} = b_i$ for $1 \leq i \leq k - 1$ (or $k=1$).
Note that $D_{t,k}$
has maximum degree $\Delta(D_{t,k})\leq 2t$, order $t^k$ and diameter $k$. 
It thus provides a nontrivial lower bound in the study of so called \emph{Moore bound}, concerning the largest order of a graph with given maximum degree and diameter, see e.g. a survey by Miller and {\v S}ir{\'a}{\v n}~\cite{Mirka}.

Now we construct our family of examples witnessing the necessity of minimum degree assumption in Conjecture~\ref{Przybylo_main_conjecture_total}
(and Theorems~\ref{first_attempt_total_theorem},~\ref{przybylo_main_total_bound}).
Given positive integers $N$ (which shall be required to be large enough
later on) and $r\geq 4$, we define $G^1_{r,N}$ (and $G^0_{r,N}$) to be the graph obtained
by taking one copy of $D_{N(r-2)^2,r-2}$, which has $[N(r-2)^2]^{r-2}$ vertices,
and $2N[N(r-2)^2]^{r-2}$ disjoint copies of $K_{r-1}$ (or $K_{r-2}$, resp.),
and identifying exactly one vertex from each of these complete graphs with some vertex of
our fixed $D_{N(r-2)^2,r-2}$ so that its every vertex is incident with exactly $2N$ 
copies of $K_{r-1}$ ($K_{r-2}$ resp.).
Note that such graph has diameter $r$, i.e., every two its vertices are $r$-neighbours,
maximum degree at most $2N(r-2)^2+2N(r-2)=2N(r-1)(r-2)$ (or $2N(r-2)^2+2N(r-3)$, resp.), and contains at least $2N[N(r-2)^2]^{r-2}(r-2)$ (or $2N[N(r-2)^2]^{r-2}(r-3)$, resp.)
vertices of degree $r-2$ ($r-3$, resp.).

\begin{observation}\label{minimium_degree_obs_total}
Let $r>3$ be an integer. Suppose $\chi''_{a,r}(G)\leq \Delta+C$ (or at least $\chi''_{a,r}(G)\leq (1+o(1))\Delta$)
for every graph $G$ with $\delta(G)\geq \delta_0$ and $\Delta(G)=\Delta$, where $C$ is some constant.
Then we must have $\delta_0\geq r-1$ if $r\geq 7$, or at least $\delta_0\geq r-2$ otherwise.
\end{observation}
\begin{pf}
Given any fixed graph and an integer $d$, let $n_d$ denote the number of vertices of degree $d$ in this graph.

Assume first that $r\geq 9$, $G$ is the graph $G^1_{r,N}$ (for some sufficiently large $N$)
and suppose we wish to construct an $r$-distant set distinguishing total colouring of $G$
using at most $\Delta+C$ (or  $(1+o(1))\Delta$) colours,
where $C$ is some constant. Note that with these many colours admitted,
there are at most ${\Delta+C \choose r-1} \leq {2N(r-1)(r-2)+C \choose r-1}$
potential pallets for vertices of degree $r-2$ in $G$.
Therefore, in order to prove that there is no desired colouring, it is sufficient to prove
that $n_{r-2}$ is larger than this quantity in our graph, i.e., that:
\begin{equation}\label{delta_lower_bound_ineq_total}
2N[N(r-2)^2]^{r-2}(r-2) > {2N(r-1)(r-2)+C \choose r-1}.
\end{equation}
Since by~(\ref{Choose_inequalities_total}),
\begin{equation}\label{delta_lower_bound_ineq_2_total}
{2N(r-1)(r-2)+C \choose r-1}\leq\left(\frac{e[2N(r-1)(r-2)+C]}{r-1}\right)^{r-1}\leq [5.5N(r-2)]^{r-1}
\end{equation}
(for $N$ sufficiently large),
by (\ref{delta_lower_bound_ineq_total}) and (\ref{delta_lower_bound_ineq_2_total}) it is then sufficient to prove the inequality
$$2N^{r-1}(r-2)^{2r-3} > [5.5N(r-2)]^{r-1},$$
or equivalently
$$\left(\frac{r-2}{5.5}\right)^{r-2} > 2.75,$$
which holds for $r\geq 9$ (the left hand side above is an increasing function of $r$ for $r\geq 9$,
as $r-2>5.5$ then and $(\frac{9-2}{5.5})^{9-2}\approx 5.41>2.75$).

Note that the same argument holds even if we admit $\Delta(1+o(1))$ instead of $\Delta+C$ colours
(as inequality (\ref{delta_lower_bound_ineq_2_total}) holds also after such substitution, i.e. after substituting
$2N(r-1)(r-2)+C$ with $2N(r-1)(r-2)(1+o(1))$ in it, where $o(1)=o_\Delta(1)$ or equivalently $o(1)=o_N(1)$ for our fixed $r$).
Moreover, by more careful estimations, one can easily show that in both cases the same is also true already for $r\geq 7$
(by proving directly that~(\ref{delta_lower_bound_ineq_total}) holds for $r=7,8$ and $N$ sufficiently large).

In the remaining cases, i.e., for $r\in\{4,5,6\}$, to show that we need the assumption $\delta(G)\geq r-2$,
it is sufficient to consider $G=G^0_{r,N}$ (instead of $G^1_{r,N}$), 
and observe that $n_{r-3}=\Omega(\Delta^{r-1})$ for such graph, while
with only $\Delta+C$ or $\Delta(1+o(1))$ colours admitted,
there would be just $O(\Delta^{r-2})$ available colour pallets for vertices of degree $r-3$ in it (for $N$ tending to infinity).
$\blacksquare$
\end{pf}

\section{Concluding Remarks}

As mentioned earlier, Theorem~\ref{przybylo_main_total_bound} implies Conjecture~\ref{Przybylo_main_conjecture_total} for regular graphs.
It suffices to substitute $\varepsilon = 1$ in the thesis of this theorem to obtain
the following.
\begin{corollary}\label{regular_graphs_corollary_total}
For every positive integer $r$, there exists $d_0$ such that:
$$\chi''_{a,r}(G)\leq \Delta(G)+8r+175$$
for every $d$-regular graph $G$ with $d\geq d_0$.
\end{corollary}
In fact this further proves that for each fixed $r$, $\chi''_{a,r}(G)\leq \Delta(G)+C$
for every regular graph, where $C$ is some constant dependent on $r$.
For large degrees it follows by Corollary~\ref{regular_graphs_corollary_total} above.
In the remaining cases, i.e., for any graph $G$ with $\Delta(G) < d_0$
one can easily prove that $\chi''_{a,r}(G)\leq C_0$, where $C_0=C_0(d_0,r)$ is some (large enough) constant.
It is e.g. sufficient to colour the edges properly first, and then greedily, one by one choose
a colour for every vertex so that it is set distinguished from its $r$-neighbours
(since the number of these is bounded by a function of $r$ and $d_0$, we shall not need more colours than some constant dependent on these two quantities), or use some other more clever approach.
Note that a similar reasoning can be applied for every other $\varepsilon$ from Theorem~\ref{przybylo_main_total_bound} (not just for $\varepsilon=1$).

Conjecture~\ref{Przybylo_main_conjecture_total} still remains open in general.
On the other hand, if we skipped the assumption $\delta(G)\geq \delta_0$ in it,
i.e., we admitted vertices of (very) small degrees in graphs considered, in particular vertices of degree $1$,
then a general upper bound for $\chi''_{a,r}(G)$ could not be smaller than $\Omega(\Delta^{r-1})$ for every $r$.
To prove that we need these many colours it is e.g. sufficient to use a construction similar to the one presented in the previous section,
but this time using always complete graphs $K_2$ instead of $K_{r-1}$
(i.e., gluing appropriate number of hanging edges to vertices of $D_{N(r-2)^2,r-2}$).
Investigating this general setting might also result in obtaining interesting new results.


\begin{thebibliography}{99}
%
\bibitem{Aigner}
M. Aigner, E. Triesch, \emph{Irregular assignments of trees and
forests}, SIAM J. Discrete Math. 3(1990), (4), 439-449.

\bibitem{Akbari}
S. Akbari, H. Bidkhori, N. Nosrati, \emph{$r$-strong edge colorings of graphs},
Discrete Math. 306 (2006) 3005--3010.

\bibitem{AlonSpencer}
N. Alon, J.H. Spencer, \emph{The Probabilistic Method}, 2nd edition, Wiley, New York,
2000.

\bibitem{BalGLS}
P.N. Balister, E. Gy\H ori, J. Lehel, R.H. Schelp, \emph{Adjacent vertex distinguishing edge-colorings},
SIAM J. Discrete Math. 21(1) (2007) 237--250.

\bibitem{Behzad}
M. Behzad, \emph{Graphs and Their Chromatic Numbers}, Ph.D. Thesis, Michigan State University (1965).

\bibitem{BonamyPrzybylo}
M. Bonamy, J. Przyby{\l}o, \emph{On the neighbour sum distinguishing index of planar graphs}, submitted.

\bibitem{ChartrandErdosOellermann}
G. Chartrand, P. Erd\H{o}s, O.R. Oellermann, \emph{How to Define an Irregular Graph},
College Math. J. 19(1) (1988) 36--42.

\bibitem{Chartrand}
G. Chartrand, M.S. Jacobson, J. Lehel, O.R. Oellermann, S. Ruiz, F. Saba,
\emph{Irregular networks},
Congr. Numer. 64 (1988) 197--210.

\bibitem{CokerJohanson}
T. Coker, K. Johannson, \emph{The adjacent vertex distinguishing total chromatic number},
Discrete Math. 312 (2012) 2741--2750.

\bibitem{FlandrinMPSW}
E. Flandrin, A. Marczyk, J. Przyby{\l}o, J-F. Sacle, M. Wo{\'z}niak,
\emph{Neighbor sum distinguishing index},
Graphs Combin. 29(5) (2013) 1329--1336.

\bibitem{Hatami}
H. Hatami, \emph{$\Delta+300$ is a bound on the adjacent vertex distinguishing edge chromatic number},
J. Combin. Theory Ser. B 95 (2005) 246--256.

\bibitem{JansonLuczakRucinski}
S. Janson, T. {\L}uczak, A. Ruci\'nski, \emph{Random Graphs}, Wiley, New York, 2000.

\bibitem{KalKarPf}
M. Kalkowski, M. Karo\'nski, F. Pfender, \emph{A new upper bound for the irregularity strength of graphs}, SIAM J. Discrete Math. 25 (2011), (3), 1319--1321.

%
%
\bibitem{KemnitzMarangio}
A. Kemnitz, M. Marangio, \emph{$d$-strong Edge Colorings of Graphs}, Graphs Combin. 30 (2014) 183--195.

\bibitem{KemnitzMarangio_total}
A. Kemnitz, M. Marangio, \emph{$d$-strong total colorings of graphs}, Disrete Math. 338 (2015) 1690--1698.

\bibitem{DistChrSurvey}
F. Kramer, H. Kramer, \emph{A survey on the distance-colouring of graphs}, Discrete Math. 308 (2008) 422-426.

%
\bibitem{Mirka}
M. Miller, J. {\v S}ir{\'a}{\v n}, \emph{Moore graphs and beyond: A survey of the degree/diameter problem},
Electron. J. Combin. (2005) 1-61, Dynamic survey DS14.

\bibitem{MocSotak}
M. Mockov\v{c}iakov\'a, R. Sot\'ak, \emph{Arbitrarily large difference between $d$-strong chromatic
index and its trivial lower bound}, Discrete Math. 313 (2013) 2000--2006.

\bibitem{MolloyReedTotal}
M. Molloy, B. Reed, \emph{A bound on the Total Chromatic Number}, Combinatorica 18 (1998) 241--280.

\bibitem{MolloyReed}
M. Molloy, B. Reed, \emph{Graph Colouring and the Probabilistic Method}, Springer, Berlin, 2002.

\bibitem{Nierhoff}
T. Nierhoff, \emph{A tight bound on the irregularity strength of
graphs}, SIAM J. Discrete Math. 13 (2000), (3), 313-323.

\bibitem{PilsniakWozniak_total}
M. Pil\'sniak, M. Wo\'zniak, \emph{On the Total-Neighbor-Distinguishing Index by Sums},
Graphs Combin. 31 (2015) 771-782.

%
\bibitem{Przybylo_distant}
J. Przyby{\l}o, \emph{Distant irregularity strength of graphs}, Discrete Math. 313(24) (2013) 2875--2880.

\bibitem{Przybylo_distant_Hatami_edge}
J. Przyby{\l}o, \emph{Distant set distinguishing edge colourings of graphs}, submitted.

\bibitem{Przybylo_CN_1}
J. Przyby{\l}o, \emph{Neighbor distinguishing edge colorings via the Combinatorial Nullstellensatz}, SIAM J. Discrete Math. 27(3) (2013) 1313--1322.

\bibitem{Przybylo_CN_2}
J. Przyby{\l}o, T-L. Wong, \emph{Neighbor distinguishing edge colorings via the Combinatorial Nullstellensatz revisited}, J. Graph Theory,
to appear (doi: 10.1002/jgt.21852).

\bibitem{Distsnt_total_Hatami_supplementary_general}
H. Qiang, M. Li, Z. Zhang, \emph{A bound on vertex distinguishing total coloring of graphs with distance constraint for recurrent event data},
Acta Math. Appl. Sin. 34 (2011) 554--559 (in Chinese).

\bibitem{Distsnt_total_Hatami_supplementary_1}
O. Sch\"afer, \emph{d-starke Totalfärbungen von Graphen} (Master thesis), TU Braunschweig, 2012.

\bibitem{Tian}
J.J. Tian, X. Liu, Z. Zhang, F. Deng, \emph{Upper bounds on the $D(\beta)$-vertex-distinguishing edge-chromatic-numbers of graphs},
LNCS 4489 (2007) 453--456.

\bibitem{Vizing2}
V. Vizing, \emph{Some Unsolved Problems in Graph Theory}, Russian Math Surveys 23 (1968) 125--141.

\bibitem{Wang_totalDelta3}
H. Wang, \emph{On the adjacent vertex distinguishing total chromatic number of the graphs with $\Delta(G) = 3$},
J. Comb. Optim. 14 (2007) 87--109.


\bibitem{WangHuang_planar}
W. Wang, D Huang, \emph{The adjacent vertex distinguishing total coloring of planar graphs}, J. Comb. Optim. 27 (2014) 379--396.

\bibitem{WangWangK_4_total_set}
W. Wang, P. Wang, \emph{On adjacent-vertex-distinguishing total coloring of $K_4$-minor free graphs},
Sci China Ser A Math 39(12) (2009) 1462--1472.


\bibitem{WangWang2}
W. Wang, Y. Wang, \emph{Adjacent vertex distinguishing total colorings of graphs with lower
average degree}, Taiwanese Journal of Mathematics 12(4) (2008) 979--990.

\bibitem{Zhang_total}
Z. Zhang, X. Chen, J. Li, B. Yao, X. Lu, J. Wang,\emph{ On adjacent-vertex-distinguishing total coloring of
graphs}, Sci China Ser A Math 48(3) (2005) 289--299.

\bibitem{ZhangDistant}
Z. Zhang, J. Li, X. Chen, et al., \emph{$D(\beta)$-vertex-distinguishing proper edge-coloring of graphs}, Acta
Math. Sinica Chin. Ser. 49 (2006) 703--708.


\bibitem{ZhangDistant_total}
Z. Zhang, J. Li, X. Chen, B. Yao, W. Wang, P. Qiu, \emph{$D(\beta)$-vertex-distinguishing total coloring of graphs},
Science China Ser. A: Mathematics 49 (2006) 1430--1440.

\bibitem{Zhang}
Z. Zhang, L. Liu, J. Wang, \emph{Adjacent strong edge coloring of graphs}, Appl. Math. Lett. 15 (2002) 623--626.

\bibitem{Distsnt_total_Hatami_supplementary_3}
Y. Zu, X. Chen, \emph{$D(5)$-vertex distinguishing total coloring of cycle}, J. Jiamusi Univ. (Nat. Science Edition) 26 (2008) 677--679 (in Chinese).

\end{thebibliography}
\end{document}